%% file: main_version.tex
\documentclass[11pt,twoside]{article}

\usepackage{xcolor}

\usepackage{fullpage}

\usepackage{epsf}
\usepackage{fancyhdr}
\usepackage{graphics}
\usepackage{graphicx}
\usepackage{psfrag}
\usepackage{comment}

\usepackage[linesnumbered,ruled]{algorithm2e}
\DontPrintSemicolon	

\usepackage{color}
\usepackage{amsthm}
\usepackage{amsfonts}
\usepackage{amsmath}
\usepackage{bm}
\usepackage{amssymb,bbm}
\usepackage{natbib}
\usepackage{algorithmic}

\usepackage{url}
\usepackage[colorlinks=True,linkcolor=magenta,citecolor=blue,urlcolor=blue,pagebackref=true,backref=true]
{hyperref}
\renewcommand*{\backref}[1]{\ifx#1\relax \else Page #1 \fi}
\renewcommand*{\backrefalt}[4]{%
    \ifcase #1 \footnotesize{(Not cited.)}%
    \or        \footnotesize{(Cited on page~#2.)}%
    \else      \footnotesize{(Cited on pages~#2.)}%
    \fi}

\usepackage{nicefrac}

\def\half{\hbox{$1\over2$}}

\usepackage{chngpage}

 \usepackage{tabularx}%

\usepackage{enumitem}
\usepackage{booktabs}
\usepackage{pbox}

\usepackage{caption,subcaption}

\usepackage{mathtools}

\usepackage{fullpage}

\input{final_macros.tex}

\begin{document}
\begin{center}

{\bf{\LARGE{Bayesian Consistency with the Supremum Metric}}}
  
\vspace*{.2in}
{\large{
\begin{tabular}{cc}
Nhat Ho$^{\diamond}$ & Stephen G. Walker$^{\diamond, \flat}$ \\
\end{tabular}
}}
\date{}
\vspace*{.2in}

\begin{tabular}{c}
Department of Statistics and Data Sciences, University of Texas at Austin$^\diamond$, \\
Department of Mathematics, University of Texas at Austin$^\flat$ \\
\end{tabular}


\vspace*{.2in}

\begin{abstract}
We present simple conditions for Bayesian consistency in the supremum metric. The key to the technique is a triangle inequality which allows us to explicitly use weak convergence, a consequence of the standard Kullback--Leibler support condition for the prior. A further condition is to ensure that smoothed versions of densities are not too far from the original density, thus dealing with densities which could track the data too closely. A key result of the paper is that we demonstrate supremum consistency using weaker conditions compared to those currently used to secure $\mathbb{L}_1$ consistency.
\end{abstract}
\end{center}

\textsl{Keywords:} Prokhorov metric, sin kernel, Fourier integral theorem, weak convergence.

\section{Introduction}

Bayesian consistency remains an open topic and has seen much progress and ideas since the seminal papers of \cite{BSW1999} and \cite{GGR1999}. A dominating sufficient, but not necessary, condition is a Kullback--Leibler support condition for the prior; 
\begin{equation}\label{klcond}
\Pi\big(D(p_0,p)<\varepsilon\big)>0
\end{equation}
for all $\varepsilon>0$. Here $D(p_0,p)=\int p_0\,\log(p_0/p)$ denotes the Kullback--Leibler divergence between $p_0$ and $p$ and $p_0$ represents the true density function from which the identically distributed $(X_i)_{i=1:n}$ are observed. Further, we write $\Pi(dp)$ to denote the prior distribution on a space of probability density functions; say $\mathbb{P}$.

It is well known that condition (\ref{klcond}) is not sufficient for strong consistency. Strong consistency holds if
\begin{equation}\label{strong}
\Pi_n(A_\varepsilon) : = \Pi(A_\varepsilon\mid X_{1:n})\to 0\quad\mbox{a.s.}\quad\,\,P_0^\infty
\end{equation}
for all $\varepsilon>0$, where $A_{\varepsilon} =\{p:d_H(p_0,p)>\varepsilon\}$ and $d_H$ is the Hellinger distance between $p_0$ and $p$. Note the Hellinger distance is equivalent to the $\mathbb{L}_1$ distance. There is a counter example in \cite{BSW1999} which shows that a posterior is not strongly consistent given only the Kullback--Leibler support condition.

The standard additional sufficient condition for consistency involves the existence of an increasing  sequence of sieves $(\mathbb{F}_n)$, which become $\mathbb{P}$ as $n\to\infty$, such that the size of $\mathbb{F}_n$, as measured by some suitable entropy, is bounded by $e^{n\kappa}$, for some $\kappa>0$, and
$\Pi(\mathbb{F}_n')<e^{-n\xi}$ for some $\xi>0$.

On the other hand, \cite{Walker2004} found a sieve, based on $\Pi$ itself, which automatically satisfies the entropy condition, and the $\mathbb{F}_n'$ condition is satisfied when 
$\sum_{j=1:\infty} \Pi(A_j)<\infty,$
where the $(A_j)_{j = 1}^{\infty}$ form a partition of $\mathbb{P}$ with respect to Hellinger neighborhoods. A recent survey of Bayesian consistency is provided in \cite{GV2017}. 

A new approach to Bayesian consistency was developed by \cite{Chae2017}. The idea is to rely on the weak convergence of the posterior and to find a minimal extension to secure strong consistency. The triangle inequality, for some strong metric $d$, yields
$$d(f_0,f)\leq d(f_0,\bar{f}_0)+d(f,\bar{f})+d(\bar{f}_0,\bar{f}),$$
where $\bar{f}$ indicates a smoothed version of $f$. Specifically in \cite{Chae2017}
$$\bar{f}(x)=\frac{1}{2h}\left[F(x+h)-F(x-h)\right]$$
is used for some smoothing parameter $h > 0$ under the univariate setting.

The triangle inequality is perfect for understanding the key aspects of strong consistency. The idea is that weak convergence can deal with the $d(\bar{f},\bar{f}_0)$ term, an assumption on $f_0$ can deal with the $d(\bar{f}_0,f_0)$ term, and a condition on $f$ not being too oscillating can deal with the $d(\bar{f},f)$ term.

In this paper we use such a strategy but alter the specific details. In particular we obtain a very simple condition for strong consistency with respect to the sup metric; i.e., $\mathbb{L}_\infty$. Previous work on the sup metric has been done by \cite{Castillo2014}, who considered contraction rates, assuming the true density on $(0,1)$ is bounded away from 0. Other papers on consistency and rates using the $\mathbb{L}_r$ metrics include~\citep{Gine2011},~\citep{Hoffmann2015} and~\citep{Scricciolo2014}.

In this work, we consider the triangle inequality
$$|f_0(x)-f(x)|\leq |f_0(x)-f_{0,R}(x)|+|f_{0,R}(x)-f_R(x)|+|f(x)-f_R(x)|$$
where $f_R$ is an alternative kernel smoothed version of $f$; specifically using the sin kernel.
That is
\begin{align}
f_R(x) : = \int_{\mathbb{R}^{d}} \prod_{j = 1}^{d} \frac{\sin(R(x_{j}-y_{j}))}{x_{j}-y_{j}}\,f(y)\,d y, \label{eq:Fourier_kernel}
\end{align}
for any $x = (x_{1}, \ldots, x_{d}) \in \mathbb{R}^{d}$.  As $R$ approaches infinity and $f \in \mathbb{L}_{1}(\mathbb{R}^{d})$, $f_R(x)$ converges to $f(x)$ according to the Fourier integral theorem~\citep{Wiener33, Bochner_1959}. 

The present paper focuses solely on consistency. The idea being that weakening the conditions on prior distributions for consistency to be achieved is and remains an important topic. These weakened conditions can then be used to achieve current rates of convergence, it is argued, with some technical applications; but the insights are coming from how the weakening of assumptions required for consistency arise.  

The layout of the paper is as follows. In Section~\ref{sec:assumptions} we outline the assumptions and initial results needed to support the general theory in Section~\ref{sec:main_results}. Section~\ref{sec:Illustration} provides an illustration of the proof strategy for establishing posterior strong consistency under widely used infinite normal mixtures. We conclude the paper with some discussion in Section~\ref{sec:discussion}.

\section{Assumptions and Initial Results}
\label{sec:assumptions}
In order to study Bayesian consistency, we define the following notion of supersmooth and ordinary smooth density functions. To simplify the presentation, $\widehat{f}$ denotes the Fourier transform of the function $f$.
\begin{definition}\label{def:generalized_tail_Fourier} 
\noindent
(1) We say that the density function $f$ is supersmooth of order $\alpha$ with scale parameter $\sigma$ if there exist universal constants $C, C_{1}$ such that for almost all $x \in \mathbb{R}^{d}$, we obtain
\begin{align*}
\abss{ \widehat{f}(x)} & \leq C \exp \parenth{ -C_{1} \sigma^2 \parenth{ \sum_{j = 1}^{d} |x_{j}|^{\alpha}} }.
\end{align*}

\noindent
(2) The density function $f$ is ordinary smooth of order $\beta$ with scale parameter $\sigma$ if there exists universal constant $c$ such that for almost all $x \in \mathbb{R}^{d}$, we have
\begin{align*}
    \abss{ \widehat{f}(x)} & \leq  c \cdot \prod_{j = 1}^{d}\frac{1}{(1 + \sigma^2 |x_{j}|^{\beta})}.
\end{align*}
\end{definition}
Examples of supersmooth functions include mixtures of location Gaussian distributions or mixture of location Cauchy distributions with similar scale parameter. In particular, when $f(x) = \sum_{i = 1}^{k} \omega_{i} \mathcal{N}(x|\mu_{i}, \sigma^2 I_{d})$ where $1 \leq k \leq \infty$, then $f$ is supersmooth density function of order 2 with scale parameter $\sigma$. When $f$ is mixture of location Cauchy distributions with scale parameter $\sigma^2 I_{d}$, then $f$ is supersmooth density function of order 1 with scale parameter $\sigma$. 

Examples of ordinary smooth functions include mixtures of location Laplace distribution with similar scale parameter $\sigma I_{d}$. In this case, these mixtures are ordinary smooth functions of order 2 with scale parameter $\sigma$.

Based on Definition~\ref{def:generalized_tail_Fourier}, we have the following result regarding the difference between $f_{R}$ and $f$. The proof of Proposition~\ref{proposition:error_Fourier_integral} is similar to that of Theorem 1 in~\cite{Ho21}; therefore, it is omitted.

\begin{proposition}
\label{proposition:error_Fourier_integral}
(1) Assume that $f$ is a supersmooth density function of order $\alpha > 0$ with scale parameter $\sigma$. Then, there exist universal constants $C$ and $C'$ such that for $R \geq C'$, we have that 
\begin{align*}
    \sup_{x \in \mathbb{R}^{d}} \abss{f_{R}(x) - f(x)}
    & \leq C \frac{R^{\max \{1 - \alpha, 0\}}}{\sigma^{2d}} \exp \parenth{-C_{1} \sigma^2 \radius^{\alpha} },
\end{align*}
where $C_{1}$ is a universal constant associated with the supersmooth density function $f$ from Definition~\ref{def:generalized_tail_Fourier}.

\noindent
(2) Assume that $f$ is a ordinary smooth density function of order $\beta > 0$ with scale parameter $\sigma$. Then, there exists a universal constants $c$ such that 
\begin{align*}
    \sup_{x \in \mathbb{R}^{d}} \abss{f_{R}(x) - f(x)}
    & \leq \frac{c}{\sigma^{2 + 2(d - 1)/\beta} R^{\beta - 1}}. 
\end{align*}
\end{proposition}

\noindent
Hence, for sufficiently large $R$, we have that 
$\sup_{x \in \mathbb{R}^{d}} |f(x)-f_R(x)|$ is sufficiently small. 

If the prior puts positive mass on all Kullback--Leibler neighborhoods of $f_0$; i.e., equation~(\ref{klcond}), then the posterior converges on weak neighborhoods of $f_0$. That is:
$$\Pi_n\left(\left|\int_{\mathbb{R}^{d}} g(y)\,(f(y)-f_0(y))\, dy\right| > \varepsilon \right)\to 0\quad\mbox{a.s.}\quad\,\,P_0^\infty$$
for all continuous and bounded functions $g$. For our purpose, we will be using the product of sinc functions, which is given by: $$g_{x,R}(y)= \prod_{j = 1}^{d} \frac{\sin R(y_{j}-x_{j})}{R(y_{j}-x_{j})},$$
which is continuous and bounded for every $x = (x_{1}, \ldots, x_{d}) \in \mathbb{R}^{d}$, $y = (y_{1}, \ldots, y_{d}) \in \mathbb{R}^{d}$ and $R > 0$.
Hence, with equation~(\ref{klcond}), we have
$$\Pi_n\left(\left|\int g_{x,R}(y)\,(f(y)-f_0(y))\, dy\right| > \varepsilon \right)\to 0\quad\mbox{a.s.}\quad\,\,P_0^\infty$$
for all $x \in \mathbb{R}^{d}$ and $R > 0$.
However, we have a stronger result to this, which is that
\begin{equation}\label{weak}
\Pi_n\left(\sup_{x \in \mathbb{R}^{d}} \left|\int g_{x,R}(y)\,(f(y)-f_0(y))\, d y\right|>\varepsilon\right)\to 0\quad\mbox{a.s.}\quad\,\,P_0^\infty
\end{equation}
for any fixed $R > 0$. This is direct to show since we can write equation~(\ref{weak}) as 
$$\Pi_n\left(\sup_{x \in \mathbb{R}^{d}} \left|\int \prod_{j=1}^d\frac{\sin(Rt_j)}{Rt_j}\,(f(t-x)-f_0(t-x))\, d t\right|>\epsilon\right)\to 0\quad\mbox{a.s.}\quad\,\,P_0^\infty.$$
This holds since if $d_P(f,f_0)\to 0$ then
$\sup_{x \in \mathbb{R}^{d}} \,d_P(f(\cdot-x),f_0(\cdot-x))\to 0,$
which follows from the definition of the Prokhorov metric $d_P$ and is tantamount to demonstrating that $(A-x)^{\epsilon}=A^{\epsilon}-x$, where 
$A^\epsilon=\{b:\|b-a\| < \epsilon,\,a\in A\}$ and  $A-x=\{a-x: a\in A\}$, which is straightforward to do.
It is worth writing this out as the following lemma:

\begin{lemma}
If $d_P(f,f_0)<\epsilon$ then $\sup_{x \in \mathbb{R}^{d}} d_P(f(\cdot-x),f_0(\cdot-x))<\epsilon$.
\end{lemma}

\vspace{0.5 em}
\noindent
\textbf{Assumption on $f_{0}$:} Throughout this paper, we utilize the following mild assumption on the true density function $f_{0}$:
\begin{equation}\label{epsr}
    \varepsilon_{R} : = \sup_{x \in \mathbb{R}^{d}} |f_{0, R}(x) - f_{0}(x)| \leq \bar{C}/R,
\end{equation}
for all $R \geq R_{0}$, where $\bar{C}$ and $R_{0}$ are some positive universal constants. This assumption is satisfied when $f_{0}$ is ordinary smooth function of order 2 with any fixed scale parameter. When $d = 1$, this condition is reduced to $f_{0}$ being almost surely second order differentiable function. 
\section{Illustration}
\label{sec:Illustration}
Before stating our general posterior consistency results in Section~\ref{sec:main_results}, we consider an illustration of our proof strategy for the posterior consistency under normal mixtures, one of the most widely used nonparametric models. To keep thing simple, we consider the normal mixture models in dimension $d=1$; whereby
$$f(x)=\sum_{j=1}^\infty w_j\,\phi((x-\mu_j)/\sigma)/\sigma,$$
the $(w_j)_{j = 1}^{\infty}$ are a set of weights, the $(\mu_j)_{j = 1}^{\infty}$ are a set of locations and the $\sigma$ is a common variance term to each normal component. Further, $\phi$ represents the usual standard normal density function. In a Bayesian model, prior distributions are assigned to the weights, locations and the variance. Our proof for the posterior consistency of normal mixtures consists of two main steps.

\vspace{0.5 em}
\noindent
\textbf{Step 1:} First, we find an appropriate upper bound for $\sup_x|f_R(x)-f(x)|$. Note that, the bound for $\sup_{x}|f_R(x)-f(x)|$ falls within the supersmooth setting in Proposition~\ref{proposition:error_Fourier_integral} and can be proved via bounding the tail of the Fourier transform of normal mixtures; nevertheless, in this section we show a different approach for deriving this bound for the normal mixture models via some closed-form computations.  

To this end, we first show that
\begin{equation}\label{ir}
I(R)=\int_{-\infty}^\infty \cos(Rx)\,\phi(x)\, dx=e^{-\half R^2}
\end{equation}
for all $R\geq 0$. Now,
$I'(R)=-\int_{-\infty}^\infty \sin(Rx)\,x\,\phi(x)\,dx,$
and using integration by parts, with $x\,\phi(x)=-\phi'(x)$, we have
$I'(R)=-R\,I(R)$ and hence equation~(\ref{ir}) holds since $I(0)=1$.

Now consider
\[
I(R)=\int_{-\infty}^\infty \cos(Rx)\,\phi(x-\mu)\,d x
=\int_{-\infty}^\infty \cos(R(x+\mu))\,\phi(x)\,d x
\]
and recall $\cos(R(x+\mu))=\cos(Rx)\cos(R\mu)-\sin(Rx)\sin(R\mu)$, 
so,
$I(R)=\cos(R\mu)\,e^{-\half R^2}$
since $\sin(Rx)$ is an odd function. Further, it is straightforward to show that
\begin{equation}\label{full}
\int_{-\infty}^\infty \cos(R(y-x))\,\phi((x-\mu)/\sigma)/\sigma\,d x=\cos(R(y-\mu))\,e^{-\half\sigma^2R^2},
\end{equation}
using suitable transforms.
If we denote
\begin{align*}
    J(R)=\int_{-\infty}^\infty \frac{\sin(Rx)}{x}\,\phi(x)\,d x,
\end{align*}
then $J'(R)$ is given by equation~(\ref{ir}), so
$J(R)=\int_{0}^{R}\,e^{-\half s^2}\,d s$
since $J(0)=0$. Hence, we find that
$$
J(y;\mu,\sigma,R)=\int_{-\infty}^\infty \frac{\sin(R(y-x))}{y-x}\,\phi((x-\mu)/\sigma)/\sigma\,d x 
=\int_0^R e^{-\half\sigma^2s^2}\,\cos(s(y-\mu))\,d s.
$$
We want to look at 
$f_R(x)-f(x)=\frac{1}{\pi}J(x;\mu,\sigma,R)-\phi((x-\mu)/\sigma)/\sigma,$
and from equation~(\ref{ir}), we have that
$$\int_0^\infty e^{-\half\sigma^2 s^2}\,\cos(s(x-\mu))\,d s=\pi\,\phi((x-\mu)/\sigma)/\sigma.$$
Therefore, for all $x \in \mathbb{R}$ we have
\[
\pi|f_R(x)-f(x)|=\left|\int_R^\infty e^{-\half\sigma^2 s^2}\,\cos(s(x-\mu))\,d s\right|\leq \int_R^\infty e^{-\half\sigma^2 s^2}\,d s
<\frac{1}{\sigma^2R}e^{-\frac{1}{2}\sigma^2R^2}.
\]
As a consequence, for any $R > 0$ we obtain that
\begin{align}
    \sup_{x \in \mathbb{R}} |f_R(x)-f(x)| < \frac{1}{\pi \sigma^2R}e^{-\frac{1}{2}\sigma^2R^2}. \label{eq:bound_normal_mixture}
\end{align}

\vspace{0.5 em}
\noindent
\textbf{Step 2:} Now, for any $R > 0$, an application of triangle inequality leads to
\begin{align*}
    \sup_{x \in \mathbb{R}} |f(x) - f_{0}(x)| \leq \sup_{x \in \mathbb{R}} |f(x) - f_{R}(x)| + \sup_{x \in \mathbb{R}} |f_{R}(x) - f_{0, R}(x)| + \sup_{x \in \mathbb{R}} |f_{0, R}(x) - f_{0}(x)|.
\end{align*}
From the assumption with $f_{0}$ in equation~\eqref{epsr}, we can rewrite the above bound as follows:
\begin{align*}
    \sup_{x \in \mathbb{R}} |f(x) - f_{0}(x)| \leq \sup_{x \in \mathbb{R}} |f(x) - f_{R}(x)| + \sup_{x \in \mathbb{R}} |f_{R}(x) - f_{0, R}(x)| + \bar{C}/R,
\end{align*}
as long as $R \geq R_{0}$. If we choose $R$ such that $R \geq \max \{2\bar{C}/\varepsilon, R_{0} \}$ then we have $\bar{C}/ R < \varepsilon/ 2$. Therefore, a direct application of union bounds shows that
\begin{align*}
    \Pi_{n}\left(\sup_{x \in \mathbb{R}} |f(x) - f_{0}(x)| > \varepsilon\right) &  \leq \Pi_{n}\left(\sup_{x \in \mathbb{R}}|f_{R}(x) - f(x)| + \sup_{x \in \mathbb{R}}|f_{R}(x) - f_{0, R}(x)| > \varepsilon/2\right)\\ 
    & \leq \Pi_{n}\left(\sup_{x \in \mathbb{R}}|f_{R}(x) - f(x)| > \varepsilon/4\right) \\
    & \hspace{8 em} + \Pi_{n}\left(\sup_{x \in \mathbb{R}}|f_{R}(x) - f_{0,R}(x)| > \varepsilon/4\right).
\end{align*}
\vspace{0.5 em}
\noindent
\textbf{Step 2.1:} For the second term $\Pi_{n}\left(\sup_{x \in \mathbb{R}}|f_{R}(x) - f_{0,R}(x)| > \varepsilon/4\right)$, equation~\eqref{weak} indicates that
\begin{align*}
    \Pi_{n}\left(\sup_{x \in \mathbb{R}} \frac{|f_{R}(x) - f_{0,R}(x)|}{R} > \varepsilon'/4\right) \to 0 \quad\mbox{a.s.}\quad\,\,P_0^\infty
\end{align*}
when $\varepsilon' = \varepsilon/ R$. It is equivalent to
\begin{align*}
    \Pi_{n}\left(\sup_{x \in \mathbb{R}} |f_{R}(x) - f_{0,R}(x)| > \varepsilon/4 \right) \to 0 \quad\mbox{a.s.}\quad\,\,P_0^\infty.
\end{align*}
\vspace{0.5 em}
\noindent
\textbf{Step 2.2:} For the first term $\Pi_{n}\left(\sup_{x \in \mathbb{R}}|f_{R}(x) - f(x)| > \varepsilon/4\right)$, the bound in equation~\eqref{eq:bound_normal_mixture} indicates that
\begin{align*}
    \Pi_{n} \left(\sup_{x \in \mathbb{R}} |f_{R}(x) - f(x)| > \varepsilon/4\right) & \leq \Pi_n \left(\frac{1}{\pi \sigma^2R}e^{-\frac{1}{2}\sigma^2R^2} > \varepsilon/4 \right) \leq \Pi_n \left(\frac{1}{\pi \sigma^2R}e^{-\frac{1}{2}\sigma^2R^2} > \frac{\bar{C}}{2R} \right),
\end{align*}
where the second inequality is due to $\varepsilon > 2 \bar{C}/ R$. Putting the above results together, for strong consistency with respect to $L_\infty$ we need to ensure that
\begin{align*}
\Pi_n \left(\frac{1}{\pi \sigma^2R^2}e^{-\frac{1}{2}\sigma^2R^2} > \frac{\bar{C}}{2R^2} \right) \to 0
\end{align*}
as $n \to \infty$. It can be achieved by taking the  prior on $\sigma^2$ to be sample size dependent and
\begin{align}
\Pi\left(\sigma^2<\tau(\tilde{C})/R^2\right)<e^{-n \cdot g(R)} \label{eq:example_normal_mixture}
\end{align}
where $\tau(\cdot)$ is the inverse function of $\exp(-z/2)/z$, $\tilde{C}$ is some universal constant, and $g(R)$ is an arbitrary choice of increasing function with $g(0)=0$ and $g(\infty)=\infty$.

The condition~\eqref{eq:example_normal_mixture} on the prior $\sigma^2$ is an instance of the general theory for strong consistency of posterior distribution for a family of supersmooth and ordinary smooth density functions. The proof strategy for this specific class of normal mixtures provides a general recipe for obtaining such general theory that we will establish in the next section. We should add that it is a common feature in the literature to consider a sequence of sample size dependent prior distributions; see for example~\citep{Ghosal-2000}.

\section{General Theory}
\label{sec:main_results}
Based on the insight from the specific class of normal mixtures in Section~\ref{sec:Illustration} and the insight on the convergence of $f_{R}$ to $f_{0,R}$ uniformly as long as $f$ converges weakly to $f_{0}$, we are now ready to state our conditions for the strong consistency with respect to the $\mathbb{L}_{\infty}$ norm. 

To simplify the presentation, we denote by $\mathcal{S}_{\alpha, \sigma}$ the set of supersmooth density function $f$ of order $\alpha > 0$ with scale parameter $\sigma > 0$. Furthermore, we define $\mathcal{O}_{\beta, \sigma}$ to be the set of ordinary smooth density function $f$ of order $\beta > 0$ with scale parameter $\sigma > 0$. We have the following result for the strong consistency with respect to $\mathbb{L}_{\infty}$ norm.

\begin{theorem}
\label{theorem:strong_consistency}
(1) (Supersmooth setting) Assume that $\alpha > 0$ and the sample size dependent prior $\Pi$ on $\sigma$ satisfies
\begin{align*}
\Pi \left(\sigma^2 \leq \tilde{C}_{1} /R^{\alpha/2} \right) \leq \exp(- \tilde{C_{2}} n \cdot g_{1}(R))
\end{align*}
for any $R \geq \tilde{R}_{1}$ where $\tilde{C}_{1}, \tilde{C}_{2}, \tilde{R}_{1}$ are some universal constants and $g_{1}$ is an increasing function in $R$ with $g_1(0)=0$ and $g_1(\infty)=\infty$. Then, we obtain that 
\begin{align*}
    \Pi_{n}\left(f \in \mathcal{S}_{\alpha, \sigma}: \sup_{x \in \mathbb{R}^{d}} \,|f(x) - f_{0}(x)| > \varepsilon \right) \to 0 \quad\mbox{a.s.}\quad\,\,P_0^\infty
\end{align*}
for all $\varepsilon > 0$.

\vspace{0.5 em}
\noindent
(2) (Ordinary smooth setting) Assume that $\beta > 1$ and the sample size dependent prior $\Pi$ on $\sigma$ satisfies that
\begin{align*}
    \Pi \left(\sigma^{2 + 2(d - 1)/ \beta} < \bar{c}_{1}/R^{(\beta - 1)/2} \right) < \exp(- \bar{c}_{2} n \cdot g_{2}(R)),
\end{align*}
for any $R \geq \bar{R}_{1}$ where $\bar{c}_{1}, \bar{c}_{2}, \bar{R}_{1}$ are some universal constants and $g_{2}$ is an increasing function in $R$ with $g_2(0)=0$ and $g_2(\infty)=\infty$. Then, we obtain that 
\begin{align*}
    \Pi_{n}\left(f \in \mathcal{O}_{\beta, \sigma}: \sup_{x \in \mathbb{R}^{d}} \,|f(x) - f_{0}(x)| > \varepsilon \right) \to 0 \quad\mbox{a.s.}\quad\,\,P_0^\infty
\end{align*}
for all $\varepsilon > 0$.
\end{theorem}

Before putting the proof we again emphasize the weakness of the conditions to secure the supremum consistency, weaker even than those currently used for the $\mathbb{L}_1$ consistency.

\begin{proof}
(1) We first consider the supersmooth setting. The proof argument is a generalization of that in Section~\ref{sec:Illustration} for normal mixtures. Here, we provide the proof for the completeness.
For any $R > 0$ and $x\in\mathbb{R}^d$, an application of triangle inequality leads to
\begin{align*}
    \sup_{x \in \mathbb{R}^{d}} |f(x) - f_{0}(x)| \leq \sup_{x \in \mathbb{R}^{d}} |f(x) - f_{R}(x)| + \sup_{x \in \mathbb{R}^{d}} |f_{R}(x) - f_{0, R}(x)| + \sup_{x \in \mathbb{R}^{d}} |f_{0, R}(x) - f_{0}(x)|.
\end{align*}
From equation~(\ref{epsr}), by choosing $R \geq \max\{2\bar{C}/\varepsilon, R_{0}\}$, we have $\varepsilon_{R} < \varepsilon/ 2$. A direct application of union bounds leads to
\begin{align*}
    \Pi_{n}\left(\sup_{x \in \mathbb{R}^{d}} |f(x) - f_{0}(x)| > \varepsilon\right) &  \leq \Pi_{n}\left(\sup_{x \in \mathbb{R}^{d}}|f_{R}(x) - f(x)| + \sup_{x \in \mathbb{R}^{d}}|f_{R}(x) - f_{0, R}(x)| > \varepsilon/2\right)\\ 
    & \leq \Pi_{n}\left(\sup_{x \in \mathbb{R}^{d}}|f_{R}(x) - f(x)| > \varepsilon/4\right) \\
    & \hspace{8 em} + \Pi_{n}\left(\sup_{x \in \mathbb{R}^{d}}|f_{R}(x) - f_{0,R}(x)| > \varepsilon/4\right).
\end{align*}
For the second term $\Pi_{n}\left(\sup_{x \in \mathbb{R}^{d}}|f_{R}(x) - f_{0,R}(x)| > \varepsilon/4\right)$, equation~(\ref{weak}) indicates that
\begin{align*}
    \Pi_{n}\left(\sup_{x \in \mathbb{R}^{d}} \frac{|f_{R}(x) - f_{0,R}(x)|}{R^{d}} > \varepsilon'/4\right) \to 0 \quad\mbox{a.s.}\quad\,\,P_0^\infty
\end{align*}
when $\varepsilon' = \varepsilon/ R^{d}$. It is equivalent to
\begin{align*}
    \Pi_{n}\left(\sup_{x \in \mathbb{R}^{d}} |f_{R}(x) - f_{0,R}(x)| > \varepsilon/4 \right) \to 0 \quad\mbox{a.s.}\quad\,\,P_0^\infty.
\end{align*}
To consider the first term $\Pi_{n}\left(\sup_{x \in \mathbb{R}^{d}}|f_{R}(x) - f(x)| > \varepsilon/4\right)$; we look at the smoothness conditions. Since $f$ is supersmooth density function, from part (1) of Proposition~\ref{proposition:error_Fourier_integral}, we obtain that
\begin{align*}
    \Pi_{n} \left(\sup_{x \in \mathbb{R}^{d}} |f_{R}(x) - f(x)| > \varepsilon/4\right) & \leq \Pi \left(C \frac{R^{\max \{1 - \alpha, 0\}}}{\sigma^{2d}} \exp \parenth{-C_{1} \sigma^2 \radius^{\alpha} } > \varepsilon/4 \right) \\
    & \leq \Pi \left(C \frac{R^{\max \{1 - \alpha, 0\}}}{\sigma^{2d}} \exp \parenth{-C_{1} \sigma^2 \radius^{\alpha} } > \bar{C}/(2R) \right),
\end{align*}
where the second inequality is due to the fact that $R \geq \max\{2\bar{C}/\varepsilon, R_{0}\}$; since the inequality $$C R^{\max \{1 - \alpha, 0\}} \exp \parenth{-C_{1} \sigma^2 \radius^{\alpha} }/\sigma^{2d} > \bar{C}/(2R)$$ implies 
  $  \lim_{R \to \infty} \sigma^2 R^{\alpha}/\log R < \infty$.
Therefore, there exist $R_{1}$ and $\bar{C}_{1}$ such that $\sigma^2 \leq \bar{C}_{1} \log R/ R^{\alpha}$ as long as $R \geq R_{1}$. Putting the above results together, we find that
\begin{align*}
\Pi \left(C \frac{R^{\max \{1 - \alpha, 0\}}}{\sigma^{2d}} \exp \parenth{-C_{1} \sigma^2 \radius^{\alpha} } > \bar{C}/(2R) \right) \leq \Pi \left(\sigma^2 \leq \bar{C}_{1} \log R/R^{\alpha} \right)
\end{align*}
as long as $R \geq \max\{2\bar{C}/\varepsilon, R_{0}, R_{1}\}$. Since $\log(R)/ R^{\alpha/2} \to 0$ as $R \to \infty$, we can find $R_{2}$ such that $\log R/ R^{\alpha/2} \leq \tilde{C}_{1}/ \bar{C}_{1}$ as long as $R \geq R_{2}$. Collecting these results, when $R \geq \max\{2\bar{C}/\varepsilon, R_{0}, R_{1}, R_{2}\}$, we have
\begin{align*}
 \Pi \left(\sigma^2 \leq \bar{C}_{1} \log R/R^{\alpha} \right) \leq \Pi \left(\sigma^2 \leq \tilde{C}_{1}/R^{\alpha/2} \right) \leq \exp(- \tilde{C_{2}} n \cdot g_{1}(R)) \leq \exp(- \tilde{C_{2}} n \cdot g_{1}(R_{0})),
\end{align*}
where the final inequality is due to the increasing property of $g_{1}$. As a consequence, we have 
\begin{align*}
    \Pi_{n}\left(\sup_{x \in \mathbb{R}^{d}} \,|f(x) - f_{0}(x)| > \varepsilon \right) \to 0 \quad\mbox{a.s.}\quad\,\,P_0^\infty
\end{align*}
for all $\varepsilon > 0$.

\vspace{0.5 em}
\noindent
(2) We now move to the ordinary smooth setting. We denote $\gamma = \min \{1, (\beta - 1)/4 \}$. From equation~\eqref{epsr}, as long as $R \geq \max \{R_{0}, 1\}$, we have
\begin{align*}
    \varepsilon_{R} = \sup_{x \in \mathbb{R}^{d}} |f_{0, R}(x) - f_{0}(x)| \leq \bar{C}/R^{\gamma}.
\end{align*}
We follow the similar argument as the supersmooth case by choosing $R \geq \max\{(2 \bar{C}/ \varepsilon)^{1/\gamma}, R_{0}, 1\}$. This choice of $R$ is to guarantee that $\varepsilon_{R} < \varepsilon/ 2$. Therefore, we also obtain that
\begin{align*}
    \Pi_{n}\left(\sup_{x \in \mathbb{R}^{d}} |f(x) - f_{0}(x)| > \varepsilon\right) & \leq \Pi_{n}\left(\sup_{x \in \mathbb{R}^{d}}|f_{R}(x) - f(x)| > \varepsilon/4\right) \\
    & \hspace{8 em} + \Pi_{n}\left(\sup_{x \in \mathbb{R}^{d}}|f_{R}(x) - f_{0,R}(x)| > \varepsilon/4\right).
\end{align*}
The second term also approaches 0 via similar argument as the supersmooth case. For the first term, the result of part (2) of Proposition~\ref{proposition:error_Fourier_integral} for ordinary smooth function indicates that
\begin{align*}
    \Pi_{n} \left(\sup_{x \in \mathbb{R}^{d}} |f_{R}(x) - f(x)| > \varepsilon/4 \right) & \leq \Pi \left(\frac{c}{\sigma^{2 + 2(d - 1)/\beta} R^{\beta - 1}} > \bar{C}/(2 R^{\gamma}) \right) \\
    & = \Pi \left(\sigma^{2 + 2(d - 1)/ \beta} < 2c/(\bar{C}R^{\beta - (1 + \gamma)}) \right).
\end{align*}
Since $2c/(\bar{C} R^{\beta - (1 + \gamma)}) < \bar{c}_{1}/R^{(\beta - 1)/2}$ where $\bar{c}_{1}$ is a constant in Theorem~\ref{theorem:strong_consistency} as long as $R^{\frac{\beta - 1}{2} - \gamma} > c/ (\bar{c}_{1} \bar{C})$, it indicates that as long as $R \geq \max\{(2 \bar{C}/ \varepsilon)^{1/\gamma}, (c/ (\bar{c}_{1} \bar{C}))^{2/(\beta - 1 - 2 \gamma)}, R_{0}, 1\}$, we have
\begin{align*}
    \Pi \left(\sigma^{2 + 2(d - 1)/ \beta} < 2c/(\bar{C}R^{\beta - (1 + \gamma)}) \right) \leq \Pi \left(\sigma^{2 + 2(d - 1)/ \beta} < \bar{c}_{1}/ R^{(\beta - 1)/2} \right) < \exp( - \bar{c}_{2} n \cdot g_{2}(1)) \to 0
\end{align*}
as $n \to \infty$ where $\bar{c}_{2}$ is the universal constant in Theorem~\ref{theorem:strong_consistency} and the final inequality is due to the increasing property of $g_{2}$. Putting the above result together, we obtain that
\begin{align*}
    \Pi_{n}\left(\sup_{x \in \mathbb{R}^{d}} \,|f(x) - f_{0}(x)| > \varepsilon \right) \to 0 \quad\mbox{a.s.}\quad\,\,P_0^\infty
\end{align*}
for all $\varepsilon > 0$.
\end{proof}

\section{Discussion}
\label{sec:discussion}
At the heart of the paper is the inequality
$$\sup_{x}|f(x)-f_0(x)|\leq \sup_{x}|f_R(x)-f(x)|+\sup_{x}|f_{0,R}(x)-f_0(x)|+\sup_{x}|f_{0,R}(x)-f_R(x)|,$$
valid for all $R>0$,
where the final term can be expressed as
$$\left|\int g_R(y)\,(f(y-x)-f_0(y-x))\,d y\right|,$$
with $g_R(y)= \prod_{j = 1}^{d} \sin(Ry_{j})/y_{j}$ for any $y = (y_{1}, \ldots, y_{d})$.
The first term is about enforcing some notion of smoothness on 
$f$ and the final term is handled by weak convergence. 

For one dimensional setting, another inequality based on the triangle inequality involves using
$$f_h(x)=\frac{1}{2h}\left[F(x+h)-F(x-h)\right],$$
as used by \cite{Chae2017}. We can now determine that
$$\sup_x|f_h(x)-f(x)|\leq \sup_{|x-y|<h} |f(x)-f(y)|$$
and so if $f$ and $f_0$ belong to a H\"older class with radius $L$ and smoothness parameter $\beta$, then
$$\sup_x|f(x)-f_0(x)|\leq d_K(f,f_0)/h+2h^\beta,$$
for any $h>0$. Here 
$$d_K(f,f_0)=\sup_x|F(x)-F_0(x)|$$
is the Kolmogorov distance where $F$ and $F_{0}$ are cumulative distribution functions of $f$ and $f_{0}$. This can be upper bounded by the 
Prokhorov metric,
$$d_K(f,f_0)\leq d_P(f,f_0)\,(1+\min\{||f||_\infty,||f_0||_\infty\}).$$
See for example~\citep{Gibbs2002}.

Hence, we should also be able to demonstrate sup norm consistency for a $\beta$ H\"older class of density once we have established weak consistency. The only condition for which we might need to construct a specific suitable prior for is the required boundedness of $||f||_\infty$.

A succinct summary of the general technique using 
$$d(f,f_0)\leq d(\bar{f},f)+d(\bar{f}_0,f_0)+d(\bar{f}_0,\bar{f}),$$
where $\bar{f}$ is a smooth version of $f$, is that the last term is handled using weak convergence, the middle term is assumed to be small, and so consistency with metric $d$ follows with some condition on $d(\bar{f},f)$. For example, a samples size dependent prior would consider $\Pi(d(\bar{f},f)>\epsilon)<e^{-n\epsilon}$ for all large $n$ and $\epsilon>0$. 


Finally, we would like to mention posterior rates of convergence. It is not a difficult task to demonstrate rates of convergence equivalent to those currently appearing in the literature, and we would anticipate with weaker conditions on the types of prior considered. However, we believe this is best reported in a future paper. The aim of the present paper is merely to point out the new technique and how elegant and useful it has turned out to be.   

\bibliographystyle{plainnat}
\bibliography{Nhat1}
\end{document}

%% file: final_macros.tex
\newcommand{\parenth}[1]{\left( #1 \right)}

\newcommand{\abss}[1]{\left| #1 \right |}

\newcommand{\radius}{\ensuremath{R}}

\newtheoremstyle{named}{}{}{\itshape}{}{\bfseries}{.}{.5em}{\thmnote{#3's }#1}
\theoremstyle{named}

\theoremstyle{plain}

\newtheorem{theorem}{Theorem}
\newtheorem{proposition}{Proposition}
\newtheorem{lemma}{Lemma}

\newtheorem{definition}{Definition}

\newlength{\widebarargwidth}
\newlength{\widebarargheight}
\newlength{\widebarargdepth}

\makeatletter
\long\def\@makecaption#1#2{
        \vskip 0.8ex
        \setbox\@tempboxa\hbox{\small {\bf #1:} #2}
        \parindent 1.5em  
        \dimen0=\hsize
        \advance\dimen0 by -3em
        \ifdim \wd\@tempboxa >\dimen0
                \hbox to \hsize{
                        \parindent 0em
                        \hfil
                        \parbox{\dimen0}{\def\baselinestretch{0.96}\small
                                {\bf #1.} #2
                                }
                        \hfil}
        \else \hbox to \hsize{\hfil \box\@tempboxa \hfil}
        \fi
        }
\makeatother

\long\def\comment#1{}
\definecolor{battleshipgrey}{rgb}{0.52, 0.52, 0.51}
\definecolor{darkgray}{rgb}{0.66, 0.66, 0.66}
\definecolor{darkgreen}{rgb}{0.0, 0.2, 0.13}
\definecolor{darkspringgreen}{rgb}{0.09, 0.45, 0.27}
\definecolor{dukeblue}{rgb}{0.0, 0.0, 0.61}
\definecolor{olivedrab7}{rgb}{0.24, 0.2, 0.12}
\definecolor{darkblue}{rgb}{0.0, 0.0, 0.55}
\definecolor{darkscarlet}{rgb}{0.34, 0.01, 0.1}
\definecolor{candyapplered}{rgb}{1.0, 0.03, 0.0}
\definecolor{ao(english)}{rgb}{0.0, 0.5, 0.0}
\definecolor{applegreen}{rgb}{0.55, 0.71, 0.0}
